\documentclass[11pt]{amsart}
\usepackage{epsfig}
\usepackage{rotate}

\newcommand{\om}{\overline m}
\newcommand{\HKR}{\mathrm{HKR}}
\newcommand{\Star}{\mathrm{Star}}
\newcommand{\poly}{\mathrm{poly}}
\newcommand{\Conf}{\mathrm{Conf}}

\newcommand{\leaves}{\mathrm{leaves}}
\renewcommand{\phi}{\varphi}

\newcommand{\F}{\mathcal F}
\newcommand{\U}{\mathcal U}

\newcommand{\DG}{{\mathrm{DG}}}
\newcommand{\DR}{{\mathrm{DR}}}

\newcommand{\ndot}{\bullet}

\def\matho#1{\mathop{\mathrm{#1}}}

\newcommand*{\ato}[2]{{\genfrac{}{}{0pt}{}{#1}{#2}}}

\DeclareMathSizes{11.1}{10}{8}{6}

\newcommand\nfrac[2]
{\dfrac{\raisebox{-1pt}{\fontsize{11.1}{10pt}\selectfont$#1$}}
       {\raisebox{ 1pt}{\fontsize{11.1}{10pt}\selectfont$#2$}}}

\newcommand{\Alt}{\matho{Alt}}

\newcommand{\C}{\mathbb C}
\newcommand{\R}{\mathbb R}

\newcommand{\D}{\mathcal D}

\newcommand{\sevafig}[3]{\begin{figure}[h]\centerline{
 \epsfig{file=#1,width=#2,angle=#3}}
\bigskip\caption{}\end{figure}}

\newtheorem*{theorem}{Theorem}

\theoremstyle{remark}

\theoremstyle{definition}
\newtheorem*{defin}{Definition}

\author{Boris Shoikhet}
\title{On the $A_\infty$-Formality conjecture}%
\date{1998}
\address{IUM, 11 Bol'shoj Vlas'evskij per.,
Moscow 121002, Russia}
\email{borya@mccme.ru}

\begin{document}
\maketitle
 \sloppy
\def\pp#1#2{\nfrac{\partial#1}{\partial#2}}

\begin{abstract}
It is proved that the associative differential graded algebra of
(polynomial) polyvector fields on a vector space (may be infinite-%
dimensional) is quasi-isomorphic to the corresponding
cohomological Hochschild complex of (polynomial) functions on
this vector space as an associative differential graded algebra.
This result is an $A_\infty$-version of the Formality conjecture
of Maxim Kontsevich [K].
\end{abstract}

\section{Configuration spaces and their compactifications}

Let $n,m$ be non-negative integers satisfying the inequality
$n+m\ge 1$. We define a configuration space $\Conf_{n,m}$ which
will play the crucial role in the construction of the
$A_\infty$-quasi-isomorphism $\F\colon
T^\ndot_{poly}\to\D^\ndot_{poly}$ in Section~2. We set
\begin{multline*}
\Conf_{n,m}=\{(p_1,\dots,p_n;q_1,\dots,q_m)\mid\\
 p_i\in\R_{<0},\
q_j\in\R_{>0}\  \text{and}\ p_1<\ldots<p_n,\ q_1>\ldots>q_m\}.
\end{multline*}
$\Conf_{n,m}$ is a smooth manifold of dimension $n+m$. The group
$G^{(1)}=\{t\mapsto at, a>0, a\in\R\}$ acts on the space
$\Conf_{n,m}$. The quotient space $C_{n,m}=\Conf_{n,m}/G^{(1)}$ is a
manifold of dimension $n+m-1$. We are going to describe a
compactification $\overline C_{n,m}$ of the space~$C_{n,m}$, which
is a \emph{manifold with corners} in the sense of Maxim Kontsevich
[K], Sect.~5.2. We defined also a manifold~$C_n$, $n\ge 2$:
$$
\Conf_n=\{(p_1,\dots,p_n)\mid p_i\in\R,\ p_1<\ldots<p_n\},
$$
$G^{(2)}=\{t\mapsto at+b,\,a,b\in\R,\, a>0\}$ and
$C_n=\Conf_n/G^{(2)}$, $\dim C_n=n-2$.

The construction of the compactification $\overline C_n$ is analogous to
[K], Sect.~5.2. The strata $C_T$ are labeled by trees $T$ and
$$
\overline C_n=\bigcup\limits_{\text{labeled trees $T$}}\prod_{v\in
V_T\setminus\{\leaves\}}C_{\Star(v)}
$$
(here $V_T$ is set of all the vertices of the tree $T$ and
$\Star(v)$ is the number of edges starting at the vertex $v$).

In the case of $C_{n,m}$ we have several possibilities:

(i) points $p_{n_1},\dots,p_n$ ($n_1\le n$) and $q_{m_1},\dots,q_m$
($m_1\le m$) are close to each other and to $0\in\R$; then the
corresponding stratum is the product $C_{n_1-1,
m_1-1}\times C_{n-n_1+1,m-m_1+1}$;

(ii) points $p_{n_1},\dots,p_{n_2}$ ($n_2-n_1\ge 1$) are close to
each other and far from $0\in\R$; then the corresponding stratum
is $C_{n_2-n_1+1}\times C_{n-n_2+n_1,m}$;
analogously in the case when points
$q_{m_1},\dots,q_{m_2}$ ($m_2-m_1\ge 1$) are close to each other;

(iii) points $p_{n_1},\dots,p_n$ ($n-n_1\ge 0$) are close to each
other and to $0\in\R$, then the corresponding stratum is
$C_{n_1-1,m}\times C_{n-n_1+1,0}$; analogously for the points
$q_{m_1},\dots,q_m$ ($m-m_1\ge 0$).

These are all the strata of codimension~$1$. 
It is easy to descibe all other strata. The strata ``of the first 
level'' are obtained when there exists several groups of the points
which are close to each other.
When we are ``looking
through a magnifying glass'' on these strata we obtain strata of
``second level,'' and so on.

\section{$A_\infty$-Formality conjecture for
$A=\C[x_1,\dots,x_d]$}

\subsection{}
Let $T^\ndot_\poly(\R^d)$ be the \emph{associative}
super-commutative algebra of polyvector fields on~$\R^d$ (with
polynomial coefficients), and let $D^\ndot_\poly(\R^d)$ be the
Hochschild cohomological complex of the algebra
$\C[x_1,\dots,x_d]$, we consider $D^\ndot_\poly(\R^d)$ as an
associative $\DG$ algebra with the product
$$
(\Theta_1\cdot\Theta_2)(f_1,\dots,f_{k+l})=\Theta_1(f_1,\dots,f_k)
\cdot\Theta_2(f_{k+1},\dots,f_{k+l}).
$$
The $A_\infty$-Formality conjecture in this case states, that
associative $\DG$ algebra $D^\ndot_\poly(\R^d)$ is quasi-isomorphic
(as an associative algebra) to its cohomology. It is well-known
result (Hochschild--Konstant--Rosenberg Theorem), that
$H^\ndot(D^\ndot_\poly(\R^d))=T^\ndot_\poly(\R^d)$, and the map
$\varphi_{\HKR}\colon T^\ndot_\poly(\R^d)\to D^\ndot_\poly(\R^d)$,
$$
\varphi_{\HKR}(\xi_1\wedge\dots\wedge\xi_k)(f_1,\dots,f_k)=
\nfrac1{k!}\Alt_{f_1,
\dots,f_k}\xi_1(f_1)\cdot\ldots\cdot\xi_k(f_k)
$$
is an quasi-isomorphism of the complexes (here
$\xi_1,\dots,\xi_k$ are vector fields and
$f_1,\dots,f_k\in\C[x_1,\dots,x_d]$).

Moreover, the inducing map $T^\ndot_\poly(\R^d)\to
H^\ndot(\D_\poly^\ndot(\R^d))$ is an isomorphism of algebras
(see, for example, [KSh], Sect.~3).

We want to construct an $A_\infty$-morphism $\F\colon
T_\poly^\ndot(\R^d)\to\D^\ndot_\poly(\R^d)$, which first
component coincides with the map~$\phi_{\HKR}$. It means, that we
want to find maps
\begin{equation}
\begin{aligned}
\F_1=\phi_\HKR&\colon
T^\ndot_\poly(\R^d)\to\D^\ndot_\poly(\R^d),\\
\F_2&\colon\otimes^2T^\ndot_\poly(\R^d)\to \D^\ndot_\poly(\R^d)[-1],\\
\F_3&\colon\otimes^3T^\ndot_\poly(\R^d)\to\D^\ndot_\poly(\R^d)[-2],\\
&\hbox to 4cm{\dotfill}
\end{aligned}
\end{equation}
such that for any $n=1,2,\dots$ and for homogeneous polyvector
fields $\gamma_1,\dots,\gamma_n$ one have:
\begin{multline}
d\F_n(\gamma_1\otimes\ldots\otimes\gamma_n)-
\sum_{\ato{k,l\ge1}{k+l=n}}\pm\F_k(\gamma_1\otimes\ldots\otimes\gamma_k)
\cdot\F_l(\gamma_{k+1},\dots,\gamma_n)-\\
-\sum_{i=1,\dots,n-1}\pm\F_{n-1}
(\gamma_1\otimes\ldots\otimes\gamma_{i-1}\otimes\gamma_i\cdot
\gamma_{i+1}\otimes\gamma_{i+2}\otimes\ldots\otimes\gamma_n)=0.
\end{multline}

\subsection{}

Let us recall the construction with graphs from~[K], Sect.~6.3.

We consider oriented graphs with $n$ vertices of the first type
and $m$ vertices of the second type, the edges start at the
vertices of the first type, and these are no loops. The vertices
of the first type are labeled by the symbols $\{1,\dots,n\}$, and
the vertices of the second type are labeled by the symbols
$\{\overline1,\dots,\om\}$.

For any vertex $k$ of the first type we denote by $\Star(k)$ the
set of edges starting at the vertex~$k$.

Any such a graph~$\Gamma$ defines a map
$$
\U_\Gamma\colon\otimes^nT^\ndot_\poly(\R^d)\to
\D^\ndot_\poly(\R^d)[1-n+l]
$$
where the number of edges of~$\Gamma$ is equal to $n+m-1+l$. This
map is defined as follows.

First of all, this map has the unique nonzero component, it is
$T^{\#\Star(1)}_\poly\otimes\ldots\otimes T^{\#\Star(n)}_\poly$.
If $\gamma_i\in T^{\#\Star(i)}_\poly(\R^d)$, we are going to
define the \emph{function}
$$
\Phi=\U_\Gamma(\gamma_1\otimes\ldots\otimes\gamma_n)
(f_1\otimes\ldots\otimes f_m).
$$

Let $E_\Gamma$ be the set of the edges of the graph~$\Gamma$.

The formula for $\Phi$ is the sum over all configurations of
indices running from 1 to $d$, labeled by $E_\Gamma$:
$$
\Phi=\sum_{I\colon E_\Gamma\to\{1,\dots,d\}}\Phi_I,
$$
where $\Phi_I$ is the product over all $n+m$ vertices of
$\Gamma$ of certain partial derivatives of functions $f_j$ and of
coefficients of $\gamma_i$.

Namely, for each vertex~$i$, $1\le i\le n$ of the first type we
associate function~$\Psi_i$ on~$\R^d$ which is a coefficient of
the polyvector field~$\gamma$:
$$
\Psi_i=\langle\gamma_i,dx^{I(e^1_i)}\otimes\ldots\otimes
dx^{I(e^{k_i}_i)}\rangle
$$
(here $k_i=\#\Star (i)$, and the edges from $\Star (i)$ are
labeled by the symbols $(e^1_i,\dots,e^{k_i}_i)$).

For each vertex $j$ of second type the associated function
$\Psi_{\overline j}$ is defined as $f_j$.

Now, at each vertex of the graph $\Gamma$ we put a function on
$\R^d$ (i.\,e. $\Psi_i$ or $\Psi_{\overline j}$).

Also, on edges of the graph $\Gamma$ there are indices $I(e)$
which label coordinates in~$\R^d$. In the next step we put into
each vertex $v$ instead of function~$\Psi_v$ its partial
derivative
$$
\left(\prod_{\ato{e\in E_\Gamma}{e=(*,v)}}\partial_{I(e)}\right)\Psi_v,
$$
and then take the product over all vertices of $\Gamma$. The
result is by definition the summand~$\Phi_i$.

\subsection{Definition}
The set $G_{n,m}$ is the set of all the oriented graphs with $n$
vertices of the first type, $m$ vertices of the second type,
$n+m-1$ edges, and all the edges start at vertices of the first
type and end at the vertices of the second type.

The following are typical pictures:
\sevafig{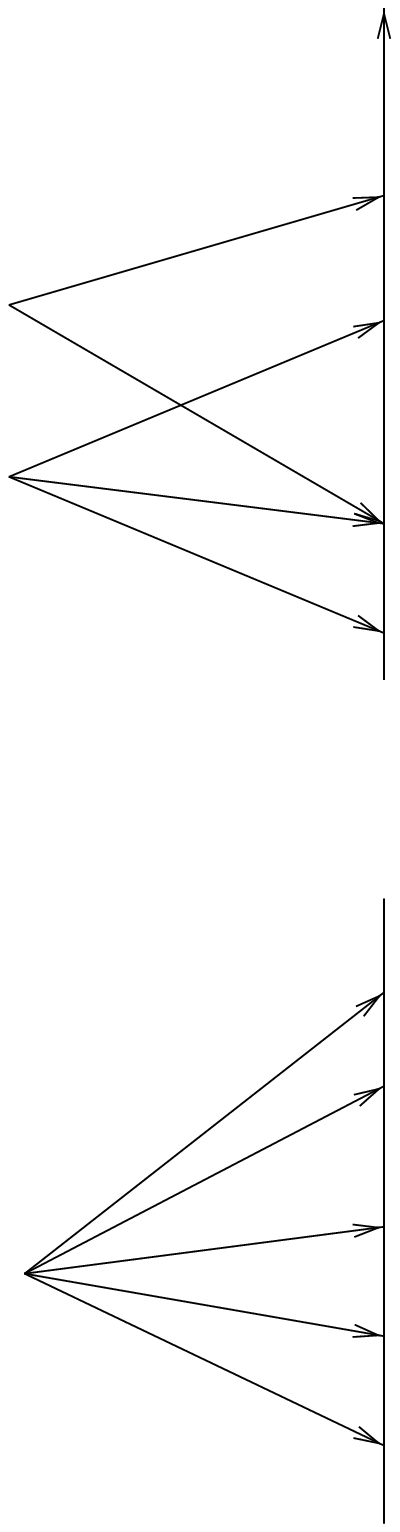}{25mm}{270}

\begin{defin}[the weight $W_\Gamma$]
Let $\phi(x)$ be any function (defined for $x<0$) such that 
is derivative~$\phi'(x)$ is a function with a compact
support, and such that $\int\phi'(x)dx=1$ (see Fig.~2).
\sevafig{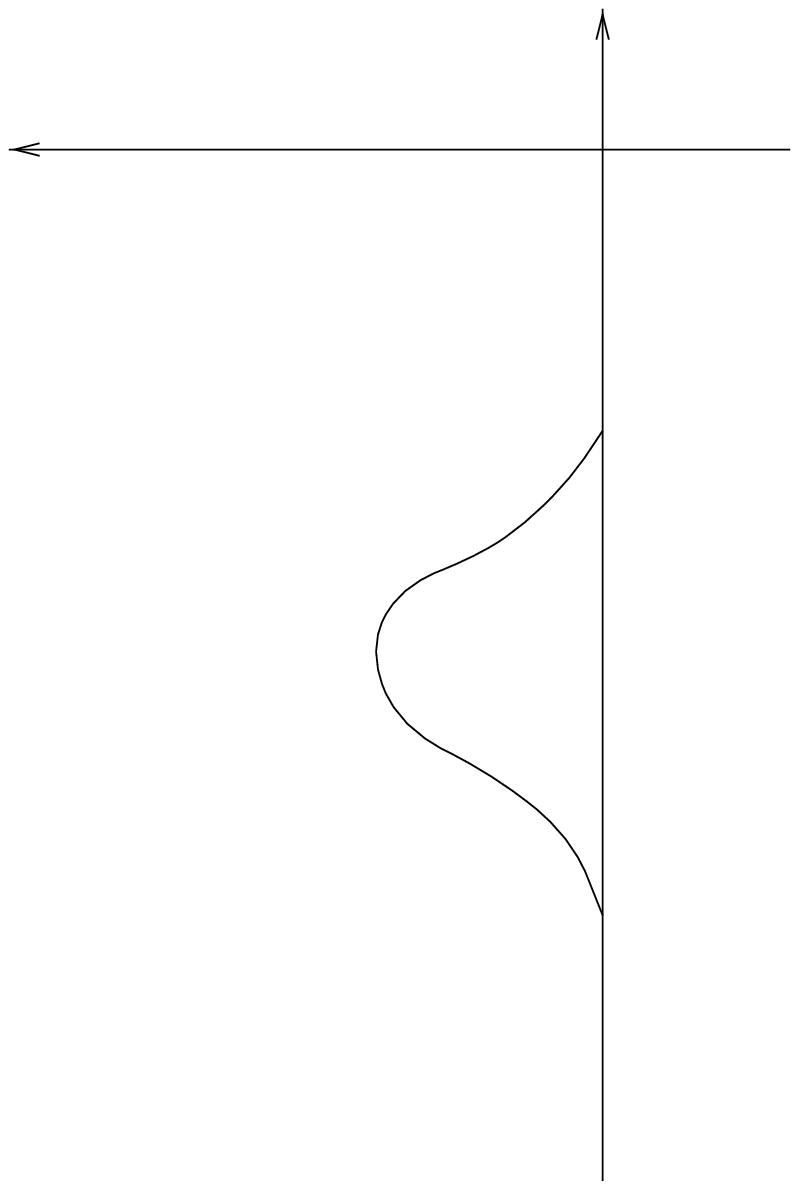}{40mm}{270}

For any pair of points $p_i,q_j$ such that $p_i<0$, $q_j>0$ we
define $\phi(p_i,q_j)$ as $\phi\left(\nfrac{p_i}{q_j}\right)$. Any
edge $e$ of $\Gamma$ defines the $1$-form $d_\DR\phi_e$ on $\overline
C_{n,m}$. By definition,
$$
W_\Gamma=\int_{\overline
C_{n,m}}\wedge_{e\in E_\Gamma}d_\DR\phi_e.
$$
\end{defin}

\begin{theorem}
The maps
$\F_1,\F_2,\F_3,\dots$ $(\F_i\colon\otimes^i T^\ndot_\poly(\R^d)\to
D^\ndot_\poly[1-i])$, where
$\F_n=\sum_{m\ge 0}\sum_{\Gamma\in G_{n,m}}W_\Gamma\times
\U_\Gamma$,
defines an $A_\infty$-morphism.
\end{theorem}

\begin{proof}

As in the $L_\infty$-case in~[K], the proof is just an
application of the Stokes formula.

Let us denote by $(F)$ the left-hand side of the formula (2). One
can write $(F)$ as a linear combination
$$
\sum_\Gamma C_\Gamma \U_\Gamma
(\gamma_1\otimes\ldots\otimes\gamma_n)(f_1\otimes\ldots\otimes f_m)
$$
where $\Gamma$ has $n$ vertices of the first type, $m$ vertices of
the second type, and $n+m-2$ edges. We want to check that
$C_\Gamma$ vanishes for each $\Gamma$. The idea (as in~[K],
Sect.~6.4) is to identify $C_\Gamma$ with the integral over the
boundary $\partial\overline C_{n,m}$ of the closed differential form
$\wedge_{e\in E_\Gamma}d\phi_e$. We have:
$$
\int_{\partial\overline C_{n,m}}\wedge_{e\in
E_\Gamma}d\phi_e=\int_{\overline C_{n,m}}d(\wedge_{e\in
E_\Gamma}d\phi_e)=0
$$

On the other hand,
$$
\int_{\partial\overline C{n,m}}\wedge_{e\in E_\Gamma}d\phi_e
$$
is an integral over all the strata of codimension $1$. All the
strata of codimension 1 were listed in the end of Section~1.

The case (i) corresponds to the second summand in~(2). The
stratum (ii) has nonzero contribution in the integral
$\int_{\partial\overline C_{n,m}}\wedge_{e\in E_\Gamma}d\phi_e$ only
in the case $\dim C_{n_2-n_1+1}=0$, i.\,e.\ it is the case where
\emph{two} neighbour points of the first type (of the second
type) are close to each other and far from 0; this case
corresponds to the third (resp., first) summand in~(2). The
stratum (iii) has a nonzero contribution in the integral
only in the case $n_1=n$ ($m_1=m$), it is the contribution
to the second (first) summand of~(2).
\end{proof}

\subsection{}
It is easy to see that $\F_1=\phi_\HKR$ (see Section~1) and
therefore the $A_\infty$-map $T^\ndot_\poly(\R^d)\to
D^\ndot_\poly(\R^d)$ we have constructed is an
$A_\infty$-quasi-isomorphism. It follows from the general theory,
that this fact is equivalent to the statement, that the $\DG$ algebras
$T^\ndot_\poly(\R^d)$ and $D^\ndot_\poly(\R^d)$ are quasi-isomorphic
in the derived category of associative $\DG$ algebras.

Note also that all the integrals $W_\Gamma$ can be easily calculated,
in the difference with the case of the $L_\infty$-Formality conjecture~[K].

\subsection{}
The remarkable difference from the case of $L_\infty$-Formality
([K]) is that our formulas make sense in the case of the algebra
of polynomials of \emph{infinite} number of variables,
$A=\C[x_1,x_2,x_3,\dots]$. The proof is the same. Also, these
formulas defines an $A_\infty$-quasi-isomorphism for any
super-algebra~$A$.

\end{document}